\documentclass[11pt,leqno,fleqn]{article}
\usepackage{epsf,amsmath,amssymb,amscd,hyperref,amsfonts}
\newcommand{\dps}{\displaystyle}

\newcommand{\T}{^{\sf T}}
\newcounter{figuren}

\newcommand{\dy}[2]{%
\refstepcounter{equation}%
\LABEL{#1}%
\begin{list}{}{
\topsep 5mm
\leftmargin 18mm
\rightmargin 0cm
\itemsep 0mm
\listparindent 0mm
\parsep 0mm
\itemsep 0mm
\labelsep 0mm
\labelwidth 18mm
}%
\item[\rm (\theequation)\hfill]
#2
\end{list}%
}
\newcommand{\dyz}[1]{%
\refstepcounter{equation}%
\begin{list}{}{
\topsep 5mm
\leftmargin 18mm
\rightmargin 0cm
\itemsep 0mm
\listparindent 0mm
\parsep 0mm
\itemsep 0mm
\labelsep 0mm
\labelwidth 18mm
}%
\item[\rm (\theequation)\hfill]
#1
\end{list}%
}
\newcommand{\dyyz}[1]{\dyz{\raggedright$\dps#1$}}

\newcommand{\de}[2]{\dy{#1}{\raggedright$\displaystyle #2 $}}
\newcommand{\dez}[1]{\dyz{\raggedright$\displaystyle #1 $}}
\newcommand{\leeg}[1]{}

\newcounter{stelling}
\newcommand{\thm}[2]{\refstepcounter{stelling}\vspace{4mm}\noindent{\bf Theorem \thestelling.}\label{#1}{\it #2}}

\newcounter{bewering}
\newcommand{\prop}[2]{\refstepcounter{bewering}\vspace{4mm}\noindent{\bf Proposition \thebewering.}\label{#1}{\it #2}}

\newcommand{\propz}[1]{\refstepcounter{bewering}\vspace{4mm}\noindent{\bf Proposition \thebewering.}{\it #1}}

\newcounter{sectie}

\newcommand{\sectz}[1]{\refstepcounter{sectie}
\section*{\boldmath \thesectie. #1}%
}
\newcounter{lit}

\newcommand{\pf}{\vspace{3mm}\noindent{\bf Proof.}\ }

\newcommand{\bx}{\hspace*{\fill} \hbox{\hskip 1pt \vrule width 4pt height 8pt depth 1.5pt \hskip 1pt}

\addvspace{4mm}}

\newcommand{\bxx}{\hspace*{\fill} \hbox{\hskip 1pt \vrule width 4pt height 8pt depth 1.5pt \hskip 1pt}}

\newcommand{\rf}[1]{{\rm (\ref{#1})}}

\newcommand{\AAA}{{\cal A}}

\newcommand{\GG}{{\cal G}}

\newcommand{\II}{{\cal I}}

\newcommand{\kint}[2]{\mbox{$\int$}}

\newcommand{\NIET}[1]{}

\newcommand{\LABEL}[1]{\label{#1}}

\newcommand{\rank}{\text{\rm rank}}

\newcommand{\oC}{{\mathbb{C}}}

\newcommand{\oZ}{{\mathbb{Z}}}

\begin{document}

\begin{center}
{\LARGE\bf Characterizing partition functions of the spin model by rank growth

}

\end{center}
\vspace{1mm}
\begin{center}
{\large
\hspace{10mm}
Alexander Schrijver\footnote{ CWI and University of Amsterdam.
Mailing address: CWI, Science Park 123, 1098 XG Amsterdam,
The Netherlands.
Email: lex@cwi.nl.}}

\end{center}

\noindent
{\small{\bf Abstract.}
We characterize which graph invariants are partition functions
of a spin model over $\oC$, in terms of the rank growth of
associated `connection matrices'.

}

\sectz{Introduction}

In this paper, all graphs are undirected and finite and may have loops and multiple edges.
An edge connecting vertices $u$ and $v$ is denoted by $uv$.
Let $\GG$ denote the collection of all undirected graphs, two of them
being the same if they are isomorphic.
A {\em graph invariant} is any function $f:\GG\to\oC$.
We consider a special class of graph invariants, namely partition functions
of spin models, defined as follows.

Let $n\in\oZ_+$.
Following de la Harpe and Jones [4], call any symmetric matrix $A\in\oC^{n\times n}$ a
{\em spin model} (over $\oC$), with $n$ {\em states}.
The {\em partition function} of $A$ is the function $p_A:\GG\to\oC$
defined for any graph $G=(V,E)$ by
\de{26jl11a}{
p_A(G):=\sum_{\kappa:V\to[n]}\prod_{uv\in E}A_{\kappa(u),\kappa(v)}.
}
Here and below, for $n\in\oZ_+$,
\dez{
[n]:=\{1,\ldots,n\}.
}
If $G$ has $k$ parallel vertices connecting $u$ and $v$, the
factor $A_{\phi(u),\phi(v)}$ occurs $k$ times in \rf{26jl11a}.

The graph invariants $p_A$ are motivated by parameters coming from
mathematical physics and from graph theory.
For instance, the Ising model corresponds to the matrix
\dez{
A=\left(\begin{array}{cc}
\exp(R/kT)&\exp(-R/kT)\\
\exp(-R/kT)&\exp(R/kT)
\end{array}
\right),
}
where $R$ is a positive constant, $k$ is the Boltzmann constant, and
$T$ is the temperature.
We refer to
[1],
[4],
and
[9]
for motivation and more examples, and
to
[3],
[5],
[6],
and
[7]
for related work and background.

In [7], partition functions
of spin models were characterized in terms of certain Moebius transforms
of graphs.
In the present paper, we characterize these graph invariants in terms
of the rank growth of associated `connection matrices'.
Rank growth of related connection matrices together with
positive semidefiniteness was considered in
Freedman, Lov\'asz, and Schrijver [3] to characterize spin functions of
real vertex models with weights on the states.

We describe the characterization.
A {\em $k$-marked graph} is a pair $(G,\mu)$ of a graph $G=(V,E)$ and a
function $\mu:[k]\to V$.
We call $i\in[k]$ a {\em mark} of vertex $\mu(i)$.
(We do not require that $\mu$ is injective, like for $k$-{\em labeled} graphs.
So a vertex may have several marks.)
Let $\GG_k$ be the collection of $k$-marked graphs.

If $(G,\mu)$ and $(H,\nu)$ are $k$-marked graphs, then
$(G,\mu)(H,\nu)$ is defined to be the graph obtained from
the disjoint union of $G$ and $H$ by identifying equally marked vertices
in $G$ and $H$.
(Another way of describing this is that we take the disjoint union
of $G$ and $H$, add edges connecting $\mu(i)$ and $\nu(i)$,
for $i=1,\ldots,k$, and finally contract each of these new edges.)

Let $f:\GG\to\oC$ and $k\in\oZ_+$.
The {\em $k$-th connection matrix} is the $\GG_k\times\GG_k$ matrix $C_{f,k}$
defined by
\dez{
(C_{f,k})_{(G,\mu),(H,\nu)}:=f((G,\mu)(H,\nu))
}
for $(G,\mu),(H,\nu)\in\GG_k$.

By $\emptyset$ we denote the graph with no vertices and edges.
We can now formulate the characterization.

\thm{26jl11b}{
Let $f:\GG\to\oC$.
Then $f=p_A$ for some symmetric $A\in\oC^{n\times n}$ and some $n\in\oZ_+$
if and only if
$f(\emptyset)=1$ and
there is a $c$ such that for each $k$: $\rank(C_{f,k})\leq c^k$.
}

\bigskip
Our proof utilizes the characterization of partition
functions of spin models given in
[7], which uses the Nullstellensatz.
One may alternatively apply the techniques
described in Freedman, Lov\'asz, and Schrijver [3].
With these techniques
one may also extend Theorem \ref{26jl11b} to more general structures like
directed graphs and hypergraphs.

A related theorem can be proved for the vertex model, where the roles
of vertices and edges are interchanged, using the characterization
given in Draisma, Gijswijt, Lov\'asz, Regts, and Schrijver [2] ---
see [8].

\sectz{Partitions}

As preliminary to the proof of Theorem \ref{26jl11b}, we give a
(most probably folklore) proposition on partitions.
A {\em partition} of a set $X$ is an (unordered) collection of pairwise
disjoint nonempty subsets of $X$ with union $X$.
The sets in $P$ are called the {\em classes} of $P$.
So $|P|$ is the number of classes of $P$.

Let $\Pi_n$ denote the collection of partitions of $[n]$.
We put $P\leq Q$ if $P$ is a refinement of $Q$, that is, if each class of $P$
is contained in some class of $Q$.
Then $(\Pi_n,\leq)$ is a lattice; we denote the join by $\vee$.

Let $Z$ be the `zeta matrix', i.e., the $\Pi_n\times\Pi_n$ matrix with
$Z_{P,Q}:=1$ if $P\leq Q$ and $Z_{P,Q}:=0$ otherwise.
Let $M:=Z^{-1}$ (the `Moebius matrix').

For $n\in\oZ_+$ and $x\in\oC$, we define the $\Pi_n\times\Pi_n$ matrix
$P_n(x)$ by
\dez{
(P_n(x))_{P,Q}:=x^{|P\vee Q|}
}
for $P,Q\in\Pi_n$.

\prop{5au11b}{
$P_n(x)$ is singular if and only if $x\in\{0,1,\ldots,n-1\}$.
}

\pf
Indeed, $MP_n(x)M\T$ is a diagonal matrix, with
\de{5au11a}{
(MP_n(x)M\T)_{P,Q}=
\delta_{P,Q}x(x-1)\cdots(x-|P|+1)
}
for $P,Q\in\Pi_n$.
Here $\delta_{P,Q}=1$ if $P=Q$ and $\delta_{P,Q}=0$ otherwise.
To prove \rf{5au11a}, we can assume $x\in\oZ_+$, as both sides are polynomials.
For $\phi:[n]\to[x]$, let $U_{\phi}$ be the partition
\dez{
U_{\phi}:=\{\phi^{-1}(i)\mid i\in[x], \phi^{-1}(i)\neq\emptyset\}.
}
Then, where $R$ and $S$ range over $\Pi_n$:
\dyyz{
(MP_n(x)M\T)_{P,Q}=
\sum_{R,S}M_{P,R}M_{Q,S}x^{|R\vee S|}
=
\sum_{R,S}M_{P,R}M_{Q,S}\sum_{\phi:[n]\to[x]\atop R\vee S\leq U_{\phi}}1
=
\sum_{R,S}M_{P,R}M_{Q,S}\sum_{\phi:[n]\to[x]\atop R,S\leq U_{\phi}}1
=
\sum_{\phi:[n]\to[x]}
\Big(\sum_{R\leq U_{\phi}}M_{P,R}\Big)
\Big(\sum_{S\leq U_{\phi}}M_{Q,S}\Big)
=
\sum_{\phi:[n]\to[x]}\delta_{P,U_{\phi}}\delta_{Q,U_{\phi}}
=
\delta_{P,Q}\sum_{\phi:[n]\to[x]}\delta_{P,U_{\phi}}
=
\delta_{P,Q}x(x-1)\cdots(x-|P|+1).
\bxx
}

\sectz{Proof of Theorem \ref{26jl11b}}

Necessity is easy, and can be seen as follows.
Let $A$ be a symmetric $n\times n$ matrix, define $f:=p_A$, and let $k\in\oZ_+$.
For any $k$-marked graph $(G,\mu)$ and any function $\lambda:[k]\to[n]$, define
\dyyz{
B_{(G,\mu),\lambda}=
\sum_{\kappa:V\to[n]\atop \kappa\circ\mu=\lambda}\prod_{uv\in E}A_{\kappa(u),\kappa(v)},
}
where $G=(V,E)$.
This defines the $\GG_k\times [n]^{[k]}$ matrix $B$, of rank at most $n^k$.
Then $C_{f,k}=BB\T$, so $C_{f,k}$ has rank at most $n^k$.
This shows necessity.

We next show sufficiency.
First observe that the conditions imply that
\de{16se12b}{
f(G\stackrel{.}{\cup}H)=f(G)f(H)
}
for all $G,H\in\GG$, where $G\stackrel{.}{\cup}H$ denotes the disjoint
union of $G$ and $H$.
This follows from the facts that the submatrix
\dez{
\left(\begin{array}{cc}
f(\emptyset)&f(G)\\
f(H)&f(G\stackrel{.}{\cup}H)
\end{array}
\right)
}
of $C_{f,0}$ has rank at most 1 and that $f(\emptyset)=1$.

By Theorem 1 in [7] it suffices to show that
for any graph $G=(V,E)$ with $V=[k]$ and $k>f(K_1)$ one has
\de{16se12c}{
\sum_{P\in\Pi_k}\mu_Pf(G/P)=0,
}
where $\mu_P:=M_{T,P}$ (with $M$ the Moebius matrix above),
where $T$ denotes the trivial partition of $[k]$ into singletons,
and where $G/P$ is the graph obtained from $G$ by merging each class of $P$
to one vertex (possibly creating several loops and multiple edges).

To prove \rf{16se12c}, from here on we fix an integer $k>f(K_1)$.
We can consider $\GG_k$ as a commutative semigroup, by maintaining the
marks in the product $(G,\mu)(H,\nu)$.
The semigroup has a unity, namely the $k$-marked graph ${\mathbf 1}_k$ with no edges
and $k$ distinct vertices marked $1,\ldots,k$.

Let $\oC\GG_k$ be the semigroup algebra of $\GG_k$.
We can extend $f$ linearly to $\oC\GG_k$.
Let $\II$ be the kernel of the matrix $C_{f,k}$,
which can be considered as a subset of $\oC\GG_k$.
Then $x\in\II$ if and only if $f(xy)=0$ for each $y\in\oC\GG_k$.
Hence $\II$ is an ideal in $\oC\GG_k$, and $\AAA:=\oC\GG_k/\II$
is a finite-dimensional commutative unital algebra with
$\dim(\AAA)=\rank(C_{f,k})$.
Moreover, as $f$ is $0$ on $\II$, $f$ has a (linear) quotient function $\hat f$ on $\AAA$.
By definition of $\II$,
for each nonzero $a\in\AAA$ there is a $b\in\AAA$ with
$\hat f(ab)\neq 0$.

As $|\Pi_n|$ grows superexponentially in $n$, there exists an $n$ such that $|\Pi_n|>c^{kn}$.
Fix an (arbitrary) bijection $s:[k]\times[n]\to[kn]$.
For each $P\in\Pi_n$ and $z\in \oC\GG_k$, let $\gamma_P(z)$ be
the following element of $\oC\GG_{kn}$.
For each $C\in P$, let $z_C$ be a copy of $z$.
For each $i\in[k]$ and $j\in [n]$ assign mark $s(i,j)$ to the vertex of $z_C$
that was marked $i$ in the original $z$,
where $C$ is the class of $P$ containing $j$.

Using \rf{16se12b}, it is direct to check that for any $P,Q\in\Pi_n$:
\de{16se12a}{
f(\gamma_P({\mathbf 1}_k)\gamma_Q(z))
=
\prod_{D\in P\vee Q}
f(z^{\text{number of classes of $Q$ contained in $D$}}).
}

\prop{17se12f}{
$\AAA$ is semisimple.
}

\pf
As $\AAA$ is commutative and finite-dimensional, it suffices to
show that any nilpotent element is zero.
To this end, suppose $a\in\AAA$ is nilpotent, with $a\neq 0$.
We can assume that $a^2=0$.
Then there is an $x\in\oC\GG_k$ with $x\not\in\II$ and $x^2\in\II$.
As $x\not\in\II$, $f(xy)\neq 0$ for some $y\in\oC\GG_k$.
Let $z:=xy$.
Then $f(z)\neq 0$ and $z^2\in\II$.
So $f(z^t)=0$ for all $t\geq 2$.
By scaling, we can assume that $f(z)=1$.

Then for any $P,Q\in\Pi_n$ we have by \rf{16se12a}
\dez{
f(\gamma_{P}({\mathbf 1}_k)\gamma_{Q}(z))=Z_{P,Q}.
}
As $Z$ is nonsingular,
this implies $\rank(C_{f,kn})\geq |\Pi_n|$, contradicting the fact that
$\rank(C_{f,kn})\leq c^{kn}<|\Pi_n|$.
\bx

Hence $\AAA\cong\oC^t$, where $t=\dim(\AAA)$.

\prop{26jl11j}{
If $a$ is a nonzero idempotent in $\AAA$, then $\hat f(a)$ is
a positive integer.
}

\pf
Let $a=z+\II$ with $z\in\oC\GG_k$.
So $f(z^t)=f(z)$ for all $t\geq 1$ (as $a^t=a$).
Then for all $P,Q\in\Pi_n$ we have by \rf{16se12a}
\dez{
f(\gamma_{P}({\mathbf 1}_k)\gamma_{Q}(z))=
(f(z))^{|P\vee Q|}.
}
As $|\Pi_n|>\rank(C_{f,kn})$, this implies that the matrix
$P_n(f(z))$ is singular.
So, by Proposition \ref{5au11b}, $f(z)\in\oZ_+$, and hence $\hat f(a)\in\oZ_+$.

Suppose finally that $a$ is a nonzero idempotent with $\hat f(a)=0$.
Then we can assume that $a$ is a minimal nonzero idempotent.
Hence $ab$ is a scalar multiple of $a$ for each $b$.
So $\hat f(ab)=0$ for each $b\in\AAA$, hence $a=0$.
\bx

For any partition $P$ of $[k]$, let $N_P$ be the $k$-marked graph
with vertex set $P$, no edges, and where mark $i\in[k]$ is given
to the element of $P$ that contains $i$.
Define the element $b$ of $\oC\GG_k$ by
\dez{
b:=\sum_{P\in\Pi_k}\mu_PN_P,
}
where, as above,$\mu_P=M_{T,P}$ for all
$P\in\Pi_k$ and $T$ is the partition
of $[k]$ consisting of singletons.

\propz{
$b$ is an idempotent in $\oC\GG_k$.
}

\pf
First note that $N_PN_Q=N_{P\vee Q}$.
Moreover, for each $R\in\Pi_k$:
\de{26jl11i}{
\sum_{P,Q\in\Pi_k\atop P\vee Q=R}\mu_P\mu_Q=\mu_R.
}
This follows from the uniqueness of $\mu$, since for each $S\in\Pi_k$ we have,
using $\mu_P=M_{T,P}$,
\dyyz{
\sum_{R\leq S}
\Big(
\sum_{P,Q\in\Pi_k\atop P\vee Q=R}\mu_P\mu_Q
\Big)
=
\sum_{P,Q\in\Pi_k\atop P\vee Q\leq S}\mu_P\mu_Q
=
\Big(\sum_{P\leq S}\mu_P\Big)^2
=
(\delta_{T,S})^2
=
\delta_{T,S}.
}
Since $MZ$ is the identity matrix, \rf{26jl11i} follows.
Hence
\dyyz{
b^2
=
\sum_{P,Q\in\Pi_k}
\mu_P\mu_QN_{P\vee Q}
=
\sum_R
\sum_{P,Q\in\Pi_k\atop P\vee Q=R}
\mu_P\mu_QN_R
=
\sum_R\mu_RN_R
=
b.
\bxx
}

Now, for any $x\in\oC$,
\dez{
\sum_{P\in\Pi_k}\mu_Px^{|P|}=x(x-1)\cdots(x-k+1)
}
(cf.\ [7] ---
it also can be derived from
\rf{5au11a}
and
\rf{26jl11i}).
Hence
\dyyz{
f(b)
=
\sum_{P\in\Pi_k}\mu_Pf(N_P)
=
\sum_{P\in\Pi_k}\mu_Pf(K_1)^{|P|}
=
f(K_1)(f(K_1)-1)\cdots(f(K_1)-k+1)
=
0.
}
The last equality follows from the facts
that
$f(K_1)$ is a nonnegative
integer and that $f(K_1)<k$.
So $f(b)=0$, and hence by Proposition \ref{26jl11j}, $b\in\II$.

Finally, to prove \rf{16se12c}, consider any graph $G$ with $k$ vertices, say
with vertex set $[k]$.
Let vertex $i\in[k]$ be marked by $i$.
Since $b\in\II$ we have $f(bG)=0$.
This is equivalent to \rf{16se12c}, and finishes the proof of Theorem \ref{26jl11b}.

\sectz{Final remark}

The condition in the theorem says that $\log(\rank(C_{f,k}))=O(k)$.
The proof shows that it can be relaxed to $\log(\rank(C_{f,k}))=o(k\log k)$,
while keeping the conditions that $\rank(C_{f,0})=1$ and $f(\emptyset)=1$.
This follows from the fact that if $\log(\rank(C_{f,k}))=o(k\log k)$,
then for each $k$ there exists an $n$ with $|\Pi_n|>\rank(C_{f,kn})$.
This is the property used in the proofs of Propositions
\ref{17se12f} and \ref{26jl11j}.

\section*{References}\label{REF}
{\small
\begin{itemize}{}{
\setlength{\labelwidth}{4mm}
\setlength{\parsep}{0mm}
\setlength{\itemsep}{1mm}
\setlength{\leftmargin}{5mm}
\setlength{\labelsep}{1mm}
}
\item[\mbox{\rm[1]}] L. Beaudin, J. Ellis-Monaghan, G. Pangborn, R. Shrock, 
A little statistical mechanics for the graph theorist,
{\em Discrete Mathematics} 310 (2010) 2037--2053.

\item[\mbox{\rm[2]}] J. Draisma, D. Gijswijt, L. Lov\'asz, G. Regts, A. Schrijver, 
Characterizing partition functions of the vertex model,
{\em Journal of Algebra} 350 (2012) 197--206.

\item[\mbox{\rm[3]}] M.H. Freedman, L. Lov\'asz, A. Schrijver, 
Reflection positivity, rank connectivity, and homomorphisms of graphs,
{\em Journal of the American Mathematical Society} 20 (2007) 37--51.

\item[\mbox{\rm[4]}] P. de la Harpe, V.F.R. Jones, 
Graph invariants related to statistical mechanical models:
examples and problems,
{\em Journal of Combinatorial Theory, Series B} 57 (1993) 207--227.

\item[\mbox{\rm[5]}] L. Lov\'asz, A. Schrijver, 
Dual graph homomorphism functions,
{\em Journal of Combinatorial Theory, Series A} 117 (2010) 216--222.

\item[\mbox{\rm[6]}] L. Lov\'asz, B. Szegedy, 
Limits of dense graph sequences,
{\em Journal of Combinatorial Theory, Series B} 96 (2006) 933--957.

\item[\mbox{\rm[7]}] A. Schrijver, 
Graph invariants in the spin model,
{\em Journal of Combinatorial Theory, Series B} 99 (2009) 502--511.  

\item[\mbox{\rm[8]}] A. Schrijver, 
Characterizing partition functions of the vertex model by rank growth,
preprint, 2012.

\item[\mbox{\rm[9]}] D.J.A. Welsh, C. Merino, 
The Potts model and the Tutte polynomial,
{\em Journal of Mathematical Physics} 41 (2000) 1127--1152.

\end{itemize}
}

\end{document}